\input amstex\documentstyle{amsppt}  
\pagewidth{12.5cm}\pageheight{19cm}\magnification\magstep1
\topmatter
\title On the cleanness of cuspidal character sheaves\endtitle
\author G. Lusztig\endauthor
\address{Department of Mathematics, M.I.T., Cambridge, MA 02139}\endaddress
\thanks{Supported in part by the National Science Foundation}\endthanks
\abstract{We prove the cleanness of cuspidal character sheaves in arbitrary characteristic in the few cases where
it was previously unknown.}\endabstract
\endtopmatter   
\document

\define\ul{\un l}

\define\bSi{\bar\Si}

\define\si{\sim}

\define\lb{\linebreak}

\define\op{\oplus}

\define\part{\partial}
\define\em{\emptyset}

\define\iy{\infty}
\define\m{\mapsto}
\define\do{\dots}

\define\sub{\subset}    

\define\T{\times}
\define\ti{\tilde}
\define\nl{\newline}
\redefine\i{^{-1}}

\define\un{\underline}

\define\ot{\otimes}
\define\bbq{\bar{\QQ}_l}

\define\ind{\text{\rm ind}}

\define\supp{\text{\rm supp}}

\define\a{\alpha}
\redefine\b{\beta}
\redefine\c{\chi}
\define\g{\gamma}
\redefine\d{\delta}

\define\p{\pi}
\define\ph{\phi}

\define\s{\sigma}

\redefine\D{\Delta}

\define\Si{\Sigma}

\define\Ph{\Phi}

\define\kk{\bold k}

\define\NN{\bold N}

\define\QQ{\bold Q}

\define\WW{\bold W}
\define\ZZ{\bold Z}

\define\cb{\Cal B}

\define\ce{\Cal E}
\define\cf{\Cal F}

\define\ch{\Cal H}

\define\cl{\Cal L}

\define\co{\Cal O}

\define\car{\Cal R}

\define\fB{\frak B}

\define\BBD{BBD}
\define\BS{BS}
\define\DL{DL}
\define\GP{GP}
\define\ICC{L1}
\define\CS{L2}
\define\WEU{L3}
\define\WEH{L4}
\define\XW{L5}
\define\OS{Os}
\define\OR{OR}
\define\SPA{Sp1}
\define\SKL{Sp2}
\define\SHF{Sh1}
\define\SH{Sh2}

\head 1. Statement of results\endhead
\subhead 1.1\endsubhead
Let $\kk$ be an algebraically closed field of characteristic exponent $p\ge1$. Let $G$ be a connected reductive 
algebraic group over $\kk$ with adjoint group $G_{ad}$.
It is known that, if $A$ is a cuspidal character sheaf on $G$, then $A=IC(\bSi,\ce)[\dim\Si]$ where $\Si$ is the
inverse image under $G@>>>G_{ad}$ of a single conjugacy class in $G_{ad}$, 
$\ce$ is an irreducible local system on $\Si$ equivariant under the conjugation $G$-action and $IC$ denotes the 
intersection cohomology complex. (For any subset $\g$ of $G$ we denote by $\bar\g$ the closure of $\g$ in $G$.) 
We say that $A$ is clean if $A|_{\bSi-\Si}=0$. This paper is concerned with the
following result.

\proclaim{Theorem 1.2} Any cuspidal character sheaf of $G$ is clean.
\endproclaim
By arguments in \cite{\CS, IV, \S17} it is enough to prove the theorem in the case where $G$ is almost simple, 
simply connected. In this case the theorem is proved in \cite{\CS,V,23.1} under the following assumption on $p$:
if $p=5$ then $G$ is not of type $E_8$; if $p=3$ then $G$ is not of type $E_7,E_8,F_4,G_2$; if $p=2$ then $G$ is 
not of type $E_6,E_7,E_8,F_4,G_2$. In the case where $p=5$ and $G$ is of type $E_8$ or $p=3$ and $G$ is of type 
$E_7,E_8,F_4,G_2$ or $p=2$ and $G$ is of type $E_6,E_7,G_2$, there are some cuspidal character sheaves on $G$ for
which the arguments of \cite{\CS, V, \S23} do not apply, but Ostrik \cite{\OS} found a simple proof for the 
cleanness of these cuspidal character sheaves. The proof of the theorem in the remaining case ($p=2$ and $G$ of 
type $E_8$ or $F_4$) is completed in 3.8; in the rest of this section we place ourselves in this case. 

Note that a portion of our proof relies on computer calculations (via the reference to \cite{\WEH} in 2.4(a)
and the references to \cite{\WEU}, \cite{\XW}).

\subhead 1.3\endsubhead
For any complex of $\bbq$-sheaves $K$ on $G$ let $\ch^iK$ be the $i$-th cohomology sheaf of $K$ and let 
${}^pH^iK$ be the $i$-th perverse cohomology sheaf of $K$. 
If $M$ is a perverse sheaf on $G$ and $A$ is a simple perverse sheaf on $G$ let $(A:M)$ be the number of times
that $A$ appears in a Jordan-H\"older series of $M$.
We write "$G$-local system" instead of "$G$-equivariant $\bbq$-local system for the conjugation action of $G$".
We set $\D=\dim G$. 

The next two properties are stated for future reference.

(a) {\it Let $A$ be a cuspidal character sheaf on $G$ and let $X$ be a noncuspidal character sheaf on $G$. Then
$H^*_c(G,A\ot X)=0$.}
\nl
(See \cite{\CS, II, 7.2}.)

(b) {\it Let $\g$ be a unipotent class in $G$ and let $\cl$ be an irreducible noncuspidal $G$-local system on 
$\g$. Then there exists a noncuspidal character sheaf $X$ of $G$ such that $\supp(X)\cap G_u\sub\bar\g$ and 
$X|_\g=\cl[d]$ for some $d\in\ZZ$.}
\nl
(See \cite{\ICC, 6.5}.)

\subhead 1.4\endsubhead
From the results already quoted we see that if $G'$ is the centralizer of a semisimple element $\ne1$ of $G$, 
then any cuspidal character sheaf of $G'$ is clean. Using \cite{\CS, II, 7.11} we then see that any cuspidal 
character sheaf of $G$ with non-unipotent support is clean. Thus it is enough to prove the cleanness of cuspidal 
character sheaves with support contained in $G_u$, the unipotent variety of $G$. For any $i\in\NN$ we denote
by $\g_i$ a distinguished unipotent class in $G$ of codimension $i$ (assuming that such class exists); note that 
$\g_i$ is unique if it exists.
According to Spaltenstein \cite{\SPA, p.336}, $\g_i$ carries an irreducible cuspidal local system precisely
when $i\in I$ where $I=\{10,20,22,40\}$ (type $E_8$) and $I=\{4,6,8,12\}$ (type $F_4$); this cuspidal local
system (necessarily of rank $1$) is unique (up to isomorphism) and denoted by $\ce_i$ except if $i=10$ (type 
$E_8$) and $i=4$ (type 
$F_4$) when there are two non-isomorphic irreducible cuspidal local systems on $\g_i$ denoted by $\ce_i,\ce'_i$.
We can then form the four admissible (see \cite{\CS, I, (7.1.10)}) complexes $A_i=IC(\bar\g_i,\ce_i)[\dim\g_i]$ 
($i\in I$) on $G$ and the admissible complex $A'_i=IC(\bar\g_i,\ce'_i)[\dim\g_i]$ (where $i=10$ for type $E_8$, 
$i=4$ for type $F_4$). According to Shoji \cite{\SH} these five 
admissible complexes are character sheaves on $G$; they are precisely the character sheaves on $G$ with support 
contained in $G_u$. From \cite{\CS, II, 7.9} we see that $A_{40}$ is clean (type $E_8$) and $A_{12}$ is clean 
(type $F_4$). According to Ostrik \cite{\OS}, $A_{10},A'_{10},A_{22}$ are clean (type $E_8$) and $A_4,A'_4,A_6$ 
are clean (type $F_4$). Moreover, from \cite{\OS} it follows that

(a) {\it if $G$ is of type $E_8$ and $i\in\ZZ$ then $\ch^i(A_{20})|_{\g_{22}}$ does not contain $\ce_{22}$ as a 
summand.}

\subhead 1.5\endsubhead
We show:

(a) {\it Let $i=20$ (type $E_8)$, $i=8$ (type $F_4$). Let $A=A_i$. Let $\g$ be a unipotent class of $G$. Then
$\op_j\ch^jA|_\g$ does not contain any irreducible noncuspidal $G$-local system as a direct summand.}
\nl
Assume that this is not true and let $\g$ be a unipotent class of minimum dimension such that $\op_j\ch^jA|_\g$ 
contains an irreducible noncuspidal $G$-local system, say $\cl$, as a direct summand. We can assume that 
$\ch^{i_0}A|_\g$ contains $\cl$ as a direct summand and that for $j>i_0$, $\ch^jA|_\g$ is a direct sum of 
irreducible cuspidal $G$-local systems on $\g$. Clearly,

(b) {\it for any unipotent class $\g'\sub\bar\g-\g$, $\op_j\ch^jA|_{\g'}$ is a direct sum of irreducible
cuspidal $G$-local systems.}
\nl
We can assume that $\g\sub\bar\g_i$; if $\g=\g_i$ the result is obvious so that we may assume that 
$\g\sub\bar\g_i-\g_i$. By 1.3(b) we can find a a noncuspidal character sheaf $X$ of $G$ such that 
$\supp(X)\cap G_u\sub\bar\g$ and $X|_\g=\cl^*[d]$ for some $d\in\ZZ$. By 1.3(a) we have $H^*_c(G,A\ot X)=0$. 
Since $\supp(A)\sub G_u$ it follows that $H^*_c(G_u,A\ot X)=0$. Since $\supp(X)\cap G_u\sub\bar\g$ it follows that
$H^*_c(\bar\g,A\ot X)=0$. 

We show that $H^*_c(\bar\g-\g,A\ot X)=0$. It is enough to show that for any unipotent class $\g'\sub\bar\g-\g$ we
have $H^*_c(\g',A\ot X)=0$. Using (b) we see that it is enough to show that for any irreducible cuspidal 
$G$-local system $\ce'$ on $\g'$ we have $H^*_c(\g',\ce'\ot X)=0$. We can find a cuspidal character sheaf $A'$ on
$G$ such that $\supp(A')=\bar\g'$, $A'|_{\g'}=\ce'[\dim\g']$. Then $A'$ must be $A_{40}$ or $A_{22}$ (for type 
$E_8$) and $A_{12}$ (for type $F_4$); in particular $A'$ is clean. Hence 

$H^*_c(\g',\ce'\ot X)=H^*_c(\bar\g',A'\ot X)=H^*_c(G,A'\ot X)$
\nl
and this is $0$ by 1.3(a).

From $H^*_c(\bar\g,A\ot X)=0$, $H^*_c(\bar\g-\g,A\ot X)=0$ we deduce that $H^*_c(\g,A\ot X)=0$ that is 
$H^*_c(\g,A\ot\cl^*)=0$. Let $\d=\dim\g$. We have $H^{2\d+i_0}_c(\g,A\ot\cl^*)=0$. We have a spectral sequence 
with $E_2^{r,s}=H^r_c(\g,\ch^s(A)\ot\cl^*)$ which converges to $H^{r+s}_c(\g,A\ot\cl^*)$. 

We show that $E_2^{r,s}=0$ if $s>i_0$. It is enough to show that $H^*_c(\g,\ce''\ot\cl^*)=0$ for any irreducible 
cuspidal $G$-local system $\ce''$ on $\g$. We can find a cuspidal character sheaf $A''$ on $G$ such that 
$\supp A''=\bar\g$, $A''|_\g=\ce''[\d]$. Since $\g\sub\bar\g_i-\g_i$ we see that $A''$ must be $A_{40}$ or 
$A_{22}$ (type $E_8$) or $A_{12}$ (type $F_4$) so that $A''$ is clean. Hence 

$H^*_c(\g,\ce''\ot\cl^*)=H^*_c(\bar\g,A''\ot X)=H^*_c(G,A''\ot X)$
\nl
and this is $0$ by 1.3(a).

We have also $E_2^{r,s}=0$ if $r>2\d$. It follows that 

$E_2^{2\d,i_0}=E_3^{2\d,i_0}=\do=E_\iy^{2\d,i_0}$.
\nl
But $E_\iy^{2\d,i_0}$ is a subquotient of $H^{2\d+i_0}(\g,A\ot\cl^*)$ hence it is zero. It follows that
$0=E_2^{2\d,i_0}=H^{2\d}_c(\g,\ch^{i_0}(A)\ot\cl^*)$. Since $\cl$ is a direct summand of $\ch^{i_0}(A)$ it 
follows that $H^{2\d}_c(\g,\cl\ot\cl^*)=0$. This is clearly a contradiction. Thus (a) is proved.

\subhead 1.6\endsubhead
We show:

(a) {\it Let $A$ be a cuspidal character sheaf on $G$ such that $\supp(A)=\bar\g$, $\g$ a unipotent class in $G$;
let $\ce$ be an irreducible $G$-local system on $\g$ such that $A|_\g=\ce[\d]$, $\d=\dim\g$. Let $Y$ be a 
noncuspidal character sheaf of G. Then $\op_j\ch^jY|_\g$ does not contain $\ce^*$ as a direct summand.}
\nl
Assume that $\op_j\ch^jY|_\g$ contains $\ce$ as a direct summand. We can find $i_0$ such that $\ch^{i_0}Y|_\g$ 
contains $\ce$ as a direct summand but $\ch^jY|_\g$ does not contain $\ce$ as a direct summand if $j>i_0$. We 
have $H^{2\d}_c(\g,\ce\ot\ch^{i_0}Y)\ne0$. By 1.3(a) we have $H^*_c(G,A\ot Y)=0$ hence $H^*_c(\bar\g,A\ot Y)=0$. 
We show that 

(b) $H^*_c(\bar\g-\g,A\ot Y)=0$. 
\nl
If $A$ is clean then (b) is obvious. Thus to prove (b) we may assume that $A=A_{20}$ (type $E_8$) and $A=A_8$ 
(type $F_4$). It is enough to show that for any unipotent class $\g'\sub\bar\g-\g$ we have $H^*_c(\g',A\ot Y)=0$.
It is enough to show that $H^*_c(\g',\ch^j(A)\ot Y)=0$ for any $j$. If $\ch^jA|_{\g'}=0$, this is obvious. Thus
we may assume that $\ch^jA|_{\g'}\ne0$. By 1.5(a), $\ch^jA|_{\g'}$ is a direct sum of (at least one) copies of 
irreducible cuspidal $G$-local systems on $\g'$. It follows that $\g'=\g_{40}$ (type $E_8$) and 
$\g'=\g_{12}$ (type $F_4$); we use that in type $E_8$ we have $\g'\ne\g_{22}$; see 1.4(a). It is then enough to 
show that $H^*_c(\g',\ce'\ot Y)=0$ where $\ce'$ is $\ce_{40}$ (type $E_8$) and $\ce'$ is $\ce_{12}$ (type $F_4$).
Let $A'=A_{40}$ (type $E_8$) and $A'=A_{12}$ (type $F_4$). Since $A'$ is clean we have 

$H^*_c(\g',\ce'\ot Y)=H^*_c(\bar\g',A'\ot Y)=H^*_c(G,A'\ot Y)$
\nl
and this is $0$ by 1.3(a). 

Using (b) and $H^*_c(\bar\g,A\ot Y)=0$ we deduce that $H^*_c(\g,A\ot Y)=0$ hence $H^*_c(\g,\ce\ot Y)=0$. Thus 
$H^{2\d+i_0}_c(\g,\ce\ot Y)=0$. We have a spectral sequence with $E_2^{r,s}=H^r_c(\g,\ce\ot\ch^sY)$ which 
converges to $H^{r+s}_c(\g,\ce\ot Y)$. We show that 

$E_2^{r,s}=0$ if $s>i_0$.
\nl
It is enough to show that $H^*_c(\g,\ce\ot\cl)=0$ for any noncuspidal irreducible $G$-local system $\cl$ on $\g$.
This follows by applying the argument in line $8$ and the ones following it in the proof of \cite{\CS, II, 7.8} 
to $\Si=\g$ (a distinguished unipotent class) and to $\cf=\ce\ot\cl$ (an irreducible $G$-local system on $\g$ not
isomorphic to $\bbq$).

We have also $E_2^{r,s}=0$ if $r>2\d$. It follows that  $E_2^{2\d,i_0}=E_3^{2\d,i_0}=\do=E_\iy^{2\d,i_0}$. But 
$E_\iy^{2\d,i_0}$ is a subquotient of $H^{2\d+i_0}(\g,\ce\ot Y)$ hence it is zero. It follows that
$0=E_2^{2\d,i_0}=H^{2\d}_c(\g,\ce\ot\ch^{i_0}Y)$. This contradicts $H^{2\d}_c(\g,\ce\ot\ch^{i_0}Y)\ne0$. This 
proves (a).

Note that a property like (a) appeared (in good characteristic) in the work of Shoji \cite{\SHF} and
Beynon-Spaltenstein \cite{\BS}.

\head 2. Preliminaries to the proof\endhead
\subhead 2.1\endsubhead
Let $\cb$ be the variety of Borel subgroups of $G$. Let $\WW$ be a set indexing the set of orbits of $G$ acting
on $\cb\T\cb$ by $g:(B,B')\m(gBg\i,gB'g\i)$. For $w\in\WW$ we write $\co_w$ for the corresponding $G$-orbit in
$\cb\T\cb$. Define $\ul:\WW@>>>\NN$ by $\ul(w)=\dim\co_w-\dim\cb$. Then $\WW$ has a natural structure of (finite)
Coxeter group with length function $\ul$ (see for example \cite{\WEU,0.2}); it is the Weyl group of $G$. 

For $w\in\WW$ let $\fB_w=\{(g,B)\in G\T\cb;(B,gBg\i)\in\co_w\}$. Define $\p_w:\fB_w@>>>G$ by $\p_w(g,B)=g$. Let 
$K_w=\p_{w!}\bbq$, a complex of sheaves on $G$. Let 

$\fB_{\le w}=\{(g,B)\in G\T\cb;(B,gBg\i)\in\cup_{y\le w}\co_y\}$,

$\fB_{<w}=\{(g,B)\in G\T\cb;(B,gBg\i)\in\cup_{y<w}\co_y\}$.
\nl
Define $\p_{\le w}:\fB_{\le w}@>>>G$ by $\p_{\le w}(g,B)=g$. Define $\p_{<w}:\fB_{<w}@>>>G$ by $\p_{<w}(g,B)=g$.
Let $K_{\le w}=\p_{\le w!}(IC(\fB_{\le w},\bbq))$, a complex of sheaves on $G$ (here $\bbq$ is viewed as a local
system on the open dense smooth subvariety $\fB_w$ of $\fB_{\le w}$).
Let $K_{<w}=\p_{<w!}(IC(\fB_{\le w},\bbq))$, a complex of sheaves on $G$.

\subhead 2.2\endsubhead
We show:

(a) {\it Let $y\in\WW$. We have ${}^pH^jK_y=0$ if $j<\D+\ul(y)$.}
\nl
We can assume that the result holds when $G$ is replaced by a Levi subgroup of a proper parabolic subgroup of $G$.
We can also assume that $G$ is semisimple. We first prove (a) for $y$ such that $y$ has minimal length in its 
conjugacy class. If $y$ is elliptic and it has minimal length in its conjugacy class in $\WW$ then, according to 
\cite{\XW, 0.3(c)}, $\p_y$ is affine and using \cite{\BBD, 4.1.1} we have ${}^pH^jK_y[\D+\ul(y)])=0$ if $j<0$
hence ${}^pH^{j+\D+\ul(y)}K_y=0$ if $j<0$ so that (a) holds for $y$. If $y$ is non-elliptic and it has 
minimal length in its conjugacy class in $\WW$ then, according to \cite{\GP, 3.2.7}, $y$ is contained in the 
subgroup $\WW'$ of $\WW$ generated by a proper subset of the set of simple reflections of $\WW$. Then $\WW'$ can 
be viewed as the Weyl group of a Levi subgroup $L$ of a proper parabolic subgroup $P$ of $G$. Define $K_{y,L}$ in
terms of $y,L$ in the same way as $K_y$ was defined in terms of $y,G$. For any $j$ we have

(b) $\ind_P^G({}^pH^jK_{y,L})={}^pH^{j+\D-\D'}K_y$.
\nl
where $\ind_P^G$ is as in \cite{\CS, 4.1} and $\D'=\dim L$. 
(This is proved along the same lines as \cite{\CS, I, 4.8(a)}.) If
$j<\D+\ul(y)$ we have $j'<\D'+\ul$ where $j'=j-\D+\D'$ hence ${}^pH^{j'}K_{y,L}=0$ so that 
$\ind_P^G({}^pH^{j'}K_{y,L})=0$ and 

$0={}^pH^{j'+\D-\D'}K_y={}^pH^jK_y$
\nl
so that (a) holds for $y$.  

We now prove (a) for any $y\in\WW$ by induction on $\ul(y)$. If $\ul(y)=0$ then $y=1$ has minimal length in its
conjugacy class and (a) holds. Now assume that $\ul(y)>0$ and that the result is known for $y'$ such that 
$\ul(y')<\ul(y)$. By \cite{\GP, 3.2.9} we can find a sequence $y=y_0,y_1,\do,y_t$ in $\WW$ such that 
$\ul(y_0)\ge\ul(y_1)\ge\do\ge\ul(y_t)$, $y_t$ has minimal length in its conjugacy class and for any $i\in[0,t-1]$
we have $y_{i+1}=s_iy_is_i$ for some simple reflection $s_i$. Since (a) is already known for $y_t$ it is enough 
to verify the following statement:

(c) {\it if $i\in[0,t-1]$ and (a) holds for $y=y_{i+1}$ then (a) holds for $y=y_i$.}
\nl 
If $\ul(y_i)=\ul(y_{i+1})$ then, by an argument similar to that in \cite{\WEU, 5.3}, we see that there exists an 
isomorphism $\fB_{y_i}@>\si>>\fB_{y_{i+1}}$ commuting with the $G$-actions and commuting with 
$\p_{y_i},\p_{y_{i+1}}$; hence $K_{y_i}=K_{y_{i+1}}$ and (c) follows in this case. Thus we can assume that 
$\ul(y_i)>\ul(y_{i+1})$ so that $\ul(y_i)=\ul(y_{i+1})+2$. We set $z=y_i,z'=y_{i+1},s=s_i$. For $(g,B)\in\fB_z$ 
we can find uniquely $B_1,B_2$ in $\cb$ such that $(B,B_1)\in\co_s$, $(B_1,B_2)\in\co_{z'}$,
$(B_2,gBg\i)\in\co_s$. Adapting an idea in \cite{\DL, \S1}, we define a partition $\fB_z=\fB_z^1\cup\fB_z^2$ by 
$$\fB_z^1=\{(g,B)\in\fB_z;B_2=gB_1g\i\}, \fB_z^2=\{(g,B)\in\fB_z;B_2\ne gB_1g\i\}.$$
Let $\p_z^1:\fB_z^1@>>>G$, $\p_z^2:\fB_z^1@>>>G$ be the restrictions of $\p_z$. Let $K_z^1=\p_{z!}^1\bbq$, 
$K_z^2=\p_{z!}^2\bbq$. It is enough to show that ${}^pH^jK_z^1=0$ and ${}^pH^jK_z^2=0$ if $j<\D+\ul(z)$. Now 
$(g,B)\m(g,B_1)$ is a morphism $\fB_z^1@>>>\fB_{z'}$, in fact an affine line bundle. It follows that 
$K^1_z=K_{z'}[-2]$. Thus ${}^pH^jK^1_z={}^pH^{j-2}K_{z'}$. This is $0$ for $j<\D+\ul(z)$ since 
$j-2<\D+\ul(z')$. Now $(g,B)\m(g,B_2)$ is a morphism $\fB_z^2@>>>\fB_{sz'}$, in fact a line bundle with the
zero-section removed. It follows that for any $j$ we have an exact sequence of perverse sheaves on $G$:
$${}^pH^{j-1}K_{sz}@>>>{}^pH^jK_z^2@>>>{}^pH^j(K_{sz}[-2]).$$
Since $\ul(sz')=\ul(z)-1$ we know that (a) holds for $sz$. If $j<\D+\ul(z)$ then $j-1<\D+\ul(sz')$ hence 
${}^pH^{j-1}K_{sz'}=0$ and ${}^pH^j(K_{sz'}[-2])={}^pH^{j-2}K_{sz'}=0$; the exact sequence above then shows that
${}^pH^jK_z^2=0$. This completes the inductive proof of (c) hence that of (a). (A somewhat similar strategy was 
employed in \cite{\OR} to prove a vanishing property for the cohomology of the varieties $X_w$ of \cite{\DL}; I 
thank X.He for pointing out the reference \cite{\OR} to me.) 

\subhead 2.3\endsubhead
We show:

(a) {\it Let $y\in\WW$ and let $A$ be a character sheaf on $G$ such that $(A:\op_j{}^pH^jK_{y'})=0$ for any 
$y'\in\WW$, $y'<y$. Then $(A:{}^pH^jK_y)=0$ for any $j\ne\D+\ul(y)$. Moreover, if $j=\D+\ul(y)$, there exists a 
(necessarily unique) subobject ${}^pH^jK_y^A$ of ${}^pH^jK_y$ such that ${}^pH^jK_y/{}^pH^jK_y^A$ is semisimple,
$A$-isotypic and $(A:{}^pH^jK_y^A)=0$.}
\nl
From our assumption we deduce (as in \cite{\CS, III, 12.7}) that $(A:\op_j{}^pH^jK_{<y})=0$. Hence the obvious 
morphism $\ph_j:{}^pH^jK_y@>>>{}^pH^jK_{\le y}$ satisfies $(A:\ker\ph_j)=0$, $(A:\text{\rm coker}\ph_j)=0$. In 
particular, $(A:{}^pH^jK_{\le y})=(A:{}^pH^jK_y)$ for any $j$. Since $\p_{\le y}$ is proper, ${}^pH^jK_{\le y}$
is semisimple, see \cite{\BBD}. Hence there is a unique direct sum decomposition of perverse sheaves
${}^pH^jK_{\le y}={}^pH^jK_{\le y,A}\op M$ such that ${}^pH^jK_{\le y,A}$ is semisimple, $A$-isotypic and 
$(A:M)=0$. Let 

$u:{}^pH^jK_{\le y,A}\op M@>>>{}^pH^jK_{\le y,A}$
\nl
be the first projection. The composition 

${}^pH^jK_y@>\ph_j>>{}^pH^jK_{\le y,A}\op M@>u>>{}^pH^jK_{\le y,A}$
\nl
is surjective (the image of $\ph_j$ contains 
${}^pH^jK_{\le y,A}$ since $(A:\text{\rm coker}\ph_j)=0$). Let ${}^pH^jK_y^A$ be the kernel of this composition. 
Then ${}^pH^jK_y/{}^pH^jK_y^A\cong{}^pH^jK_{\le y,A}$ hence ${}^pH^jK_y/{}^pH^jK_y^A$ is semisimple, 
$A$-isotypic. Moreover
$$\align&(A:{}^pH^jK_y^A)=(A:{}^pH^jK_y)-(A:{}^pH^jK_y^A{}^pH^jK_y)\\&
=(A:{}^pH^jK_{\le y})-(A:{}^pH^jK_{\le y,A})=(A:M)=0.\endalign$$
By the Lefschetz hard theorem \cite{\BBD, 5.4.10} we have for any $j'$:

${}^pH^{-j'}(K_{\le y}[\D+\ul(y)])\cong{}^pH^{j'}(K_{\le y}[\D+\ul(y)])$ 
\nl
hence for any $j$, ${}^pH^jK_{\le y}\cong{}^pH^{2\D+2\ul(y)-j}K_{\le y}$. It follows that

${}^pH^jK_{\le y,A}\cong{}^pH^{2\D+2\ul(y)-j}K_{\le y,A}$
\nl
so that 

(b) ${}^pH^jK_y/{}^pH^jK_y^A\cong{}^pH^{2\D+2\ul(y)-j}K_y/{}^pH^{2\D+2\ul(y)-j}K_y^A$.
\nl
\nl
Using 2.2(a) we have ${}^pH^jK_y=0$ if $j<\D+\ul(y)$. Hence ${}^pH^jK_y/{}^pH^jK_y^A=0$ if $j<\D+\ul(y)$. Using 
(b) we deduce ${}^pH^jK_y/{}^pH^jK_y^A=0$ if $j>\D+\ul(y)$. Thus ${}^pH^jK_y/{}^pH^jK_y^A=0$ if $j\ne\D+\ul(y)$. 
Since $(A:{}^pH^jK_y^A)=0$ it follows that $(A:{}^pH^jK_y)=0$ if $j\ne\D+\ul(y)$. This completes the proof of (a).

\subhead 2.4\endsubhead
In this subsection we assume that $G$ is adjoint. Let $w$ be an elliptic element of $\WW$ which has minimal 
length in its conjugacy class $C$. We assume that the unipotent class $\g=\Ph(C)$ in $G$ ($\Ph$ as in 
\cite{\WEU, 4.1}) is distinguished and that $\det(1-w)$ is a power of $p$ (the determinant is taken in the 
reflection representation of $\WW$). According to \cite{\WEH, 0.2}, 

(a) {\it the variety $\p_w\i(\g)$ is a single $G$-orbit for the $G$-action \lb
$x:(g,B)\m(xgx\i,xBx\i)$ on $\fB_w$.}
\nl
We show:

(b) {\it $K_w[2\ul(w)]|_\g\cong\op_\ce\ce^{\op\text{\rm rk}(\ce)}$ where $\ce$ runs over all irreducible 
$G$-local systems on $\g$ (up to isomorphism).}
\nl
Let $(g,B)\in\p_w\i(\g)$ and let $Z_G(g)$ be the centralizer of $g$. According to \cite{\WEU, 4.4(b)} we have 

(c) $\dim Z_G(g)=\ul(w)$.
\nl
We have a commutative diagram
$$\CD
G@>\b>>G/Z_G(g)\\
@V\a VV        @V\a'VV\\
\p_w\i(\g)@>\s>>\g
\endCD$$
where $\b$ is the obvious map, $\a(x)=(xgx\i,xBx\i)$, $\a'(x)=xgx\i$, $\s(g',B')=g'$. Now $\a$ is surjective by 
(a); it is also injective since by \cite{\WEU, 5.2} the isotropy groups of the $G$-action on $\fB_w$ are trivial 
(we use our assumption on $\det(1-w)$). Thus $\a$ is a bijective morphism so that $\a_!\bbq=\bbq$. Hence 
$$K_w[2\ul(w)]|_\g=\s_!\bbq[2\ul(w)]=\s_!\a_!\bbq[2\ul(w)]=\a'_!\b_!\bbq[2\ul(w)].$$
We now factorize $\b$ as follows: 

$G@>\b_1>>G/Z_G(g)^0@>\b_2>>G/Z_G(g)$.
\nl
Since all fibres of $\b_1$ are 
isomorphic to $Z_G(g)^0$ (an affine space of dimension $\ul(w)$, see (c)), we have $\b_!\bbq\cong\bbq[-2\ul(w)]$.
Thus
$$K_w[2\ul(w)]|_\g\cong\a'_!\b_{2!}\bbq=(\a'\b_2)_!\bbq.$$
Now $\a'\b_2$ is a principal covering with (finite) group $Z_G(g)/Z_G(g)^0$; (b) follows. (The proof above has 
some resemblance to the proof of \cite{\CS, IV, 21.11}.)

\head 3. Completion of the proof\endhead
\subhead 3.1\endsubhead
In this section (except in 3.10) we assume that $p=2$ and that $G$ is of type $E_8$ or $F_4$. We have the 
following result:

(a) {\it Let $y\in\WW$ be an elliptic element of minimal length in its conjugacy class and let $i\in I$ be such 
that $\p_y\i(\g_i)\ne\em$. Then $\ul(y)\ge i$.}
\nl
Indeed, from \cite{\WEU, 5.7(iii)} we have $\dim\g_i\ge\D-\ul(y)$ and it remains to use that $\dim\g_i=\D-i$.

\subhead 3.2\endsubhead
Let $i\ne i'$ in $I$. Then 

(a) {\it $\op_j\ch^jA_{i'}|_{\g_i}$ does not contain $\ce_i$ as a direct summand except possibly when 
$i=40,i'=20$ (type $E_8$) and $i=12,i'=8$ (type $F_4$).}
\nl
If $i'\ne20$ (type $E_8$) and $i'\ne8$ (type $F_4$) this follows from the cleanness of $A_{i'}$. If 
$i'=20,i\ne40$ (type $E_8$) and $i'=8,i\ne12$ (type $F_4$) this follows from the fact that
$\g_i\not\sub\bar\g_{i'}$ except when $i'=20,i=22$ (type $E_8$) when the result follows from 1.4(a).

Note that if $i'=10$ (type $E_8$) and $i'=4$ (type $F_4$) then

(b) $\op_j\ch^jA'_{i'}|_{\g_i}=0$
\nl
by the cleanness of $A'_{i'}$. If $i=10$ (type $E_8$) and $i=4$ (type $F_4$) then

(c) $\op_j\ch^jA_{i'}|_{\g_i}=0$
\nl
since $\g_i\not\sub\bar\g_{i'}$.

\subhead 3.3\endsubhead
We show:

(a) {\it Assume that $i\in I$, $y\in\WW$, $\ul(y)<i$. Assume also that $i\ne40$ (type $E_8$) and $i\ne12$ (type 
$F_4$). Then $(A_i:\op_j{}^pH^jK_y)=0$. If, in addition, $i=10$ for type $E_8$ and $i=4$ for type $F_4$ then 
$(A'_i:\op_j{}^pH^jK_y)=0$.}
\nl
Assume that the first assertion of (a) is false. Then we can find $y'\in\WW$ such that $\ul(y')<i$, 
$(A_i:\op_j{}^pH^jK_{y'})\ne0$ and $(A_i:\op_j{}^pH^jK_{y''})=0$ for any $y''\in\WW$ with $y''<y'$. Using 2.3(a) 
we see that $(A_i:{}^pH^jK_{y'})=0$ for any $j\ne\D+\ul(y')$ hence $(A_i:{}^pH^jK_{y'})\ne0$ for 
$j=\D+\ul(y')$. It follows that $\sum_j(-1)^j(A_i:{}^pH^jK_{y'})\ne0$. Using \cite{\CS, I, 6.5} we deduce that
$\sum_j(-1)^j(A_i:{}^pH^jK_{y'_1})\ne0$ for any $y'_1\in\WW$ that is conjugate to $y'$. If $y'$ is not elliptic
then some $y'_1$ in the conjugacy class of $y'$ is contained in the subgroup $\WW'$ of $\WW$ generated by a 
proper subset of the set of simple reflections of $\WW$. Then $\WW'$ can be viewed as the Weyl group of a Levi 
subgroup $L$ of a proper parabolic subgroup $P$ of $G$. Define $K_{y'_1,L}$ in terms of $y'_1,L$ in the same way 
as $K_y$ was defined in terms of $y,G$. We have $(A_i:\op_j{}^pH^jK_{y'_1})\ne0$. From this and from the
equality 2.2(b) (for $y'_1$ instead of $y$) we deduce that  $(A_i:\ind_P^G({}^pH^jK_{y'_1,L}))\ne0$ for some $j$.
Hence $(A_i:\ind_P^G(\ti A))\ne0$ for some character sheaf $\ti A$ on $L$; this contradicts the fact that $A_i$ 
is a cuspidal character sheaf. We see that $y'$ is elliptic. If the conjugacy class of $y'$ contains an element
$y'_2$ such that $\ul(y'_2)<\ul(y')$ then using again \cite{\CS, I, 6.5}, we deduce from 
$\sum_j(-1)^j(A_i:{}^pH^jK_{y'})\ne0$ that $\sum_j(-1)^j(A_i:{}^pH^jK_{y'_2})\ne0$ hence 
$(A_i:{}^pH^jK_{y'_2})\ne0$, contradicting the choice of $y'$. We see that $y'$ has minimal length in its 
conjugacy class.

For any $G$-equivariant perverse sheaf $M$ on $G$ we set $\c_i(M)=\sum_j(-1)^j(\ce_i:\ch^jM|_{\g_i})$ where 
$(\ce_i:?)$ denotes multiplicity in a $G$-local system. For any noncuspidal character sheaf $X$ on $G$ we have 
$\c_i(X)=0$, see 1.6(a). For any cuspidal character sheaf $X$ on $G$ with nonunipotent support we have clearly 
$\c_i(X)=0$. 

If $i'\in I-\{i\}$ then $\c_i(A_{i'})=0$ by 3.2. Also, if $i'=10$ (type $E_8$) and $i'=4$ (type $F_4$) and 
$i'\ne i$ then $\c_i(A'_{i'})=0$ by 3.2. 

From the definition we have $\c_i(A_i)\ne0$. Since $(A_i:{}^pH^jK_{y'})=0$ for any $j\ne\D+\ul(y')$ and 
$(A_i:{}^pH^jK_{y'})\ne0$ for $j=\D+\ul(y')$ it follows that $\c_i({}^pH^jK_{y'})\ne0$ for $j\ne\D+\ul(y')$ and 
$\c_i({}^pH^jK_{y'})=0$ for $j=\D+\ul(y')$. Hence $\sum_j(-1)^j\c_i({}^pH^jK_{y'})\ne0$. Hence 
$\sum_j(-1)^j(\ce_i:\ch^jK_{y'}|_{\g_i})\ne0$. It follows that $K_{y'}|_{\g_i}\ne0$ so that
$\p_{y'}\i(\g_i)\ne\em$. Using 3.1(a) we deduce that $\ul(y')\ge i$. This contradicts $\ul(y')<i$ and proves the 
first assertion of (a). The proof of the second assertion of (a) is entirely similar,

\subhead 3.4\endsubhead
We now prove a weaker version of 3.3(a) assuming that $i=40$ (type $E_8$) and $i=12$ (type $F_4$).

(a) {\it If $y\in\WW$, $\ul(y)<20$ (type $E_8$) and $\ul(y)<8$ (type $F_4$) then $(A_i:\op_j{}^pH^jK_y)=0$.}
\nl
We go through the proof of 3.3(a). The first two paragraphs remain unchanged. In the third pragraph, the sentence 

"If $i'\in I-\{i\}$ then $\c_i(A_{i'})=0$ by 3.2." 
\nl
must be modified as follows: 

"If $i'\in I-\{i\}$ and $i'\ne20$ (type $E_8$) and $i'\ne8$ (type $F_4$) then
$\c_i(A_{i'})=0$ by 3.2. Moreover, if $i'=20$ (type $E_8$) and $i'=8$ (type $F_4$) then by 3.3(a),
$(A_{i'}:{}^pH^jK_{y'})=0$ for any $j$, since $\ul(y')<20$ (type $E_8$) and $\ul(y')<8$ (type $F_4$)".
\nl
Then the fourth paragraph remains unchanged and (a) is proved.

\subhead 3.5\endsubhead
For any $i\in I$ we consider the conjugacy class $C_i$ of $\WW$ whose elements have the following characteristic 
polynomial in the reflection representation $\car$ of $\WW$: 

(type $E_8$): $q^8-q^4+1$ (if $i=10$), $(q^4-q^2+1)^2$ (if $i=20$), $(q^2-q+1)^2(q^4-q^2+1)$ (if $i=22$),
$(q^2-q+1)^4$ (if $i=40$);

(type $F_4$): $(q^4-q^2+1)$ (if $i=4$), $q^4+1$ (if $i=6$), $(q^2-q+1)^2$ (if $i=8$), $(q^2+1)^2$ (if $i=12$).
\nl
We choose an element $w_i$ of minimal length in $C_i$. Then $\ul(w_i)=i$. Note that $w_i$ is elliptic and 
$\det(1-w_i,\car)$ is $1$ (type $E_8$) and a power of $2$ (type $F_4$). 

Let $\Ph$ be the (injective) map from the set of elliptic conjugacy classes in $\WW$ to the set of unipotent 
classes in $G$ defined in \cite{\WEU, 4.1}. We have $\Ph(C_i)=\g_i$.

Note that the correspondence between $C_i$ and (the characteristic zero analogue of) $\g_i$ appeared in another 
context in the (partly conjectural) tables of Spaltenstein \cite{\SKL}.

\subhead 3.6\endsubhead
In this subsection we set $i=20,i'=40$ (type $E_8$) and $i=8,i'=12$ (type $F_4$). Let $w=w_i$. We have the 
following results.

(a) {\it If $j\ne\D+i$ then $(A_i:{}^pH^jK_w)=0$ and $(A_{i'}:{}^pH^jK_w)=0$.}

(b) {\it If $j=\D+i$ then $(A_i:{}^pH^jK_w)=1$; there exists a unique subobject $Z$ of ${}^pH^jK_w$ such that
$(A_i:Z)=0$ and ${}^pH^jK_w/Z\cong A_i$ and there exists a unique subobject $Z'$ of ${}^pH^jK_w$ such that 
$(A_{i'}:Z')=0$ and ${}^pH^jK_w/Z'$ is semisimple, $A_{i'}$-isotypic.}
\nl
(a) follows from 2.3(a) applied with $y=w$ and with $A$ equal to $A_i$ or $A_{i'}$. (The assumptions of 2.3(a) 
are satisfied by 3.3(a), 3.4(a).) 
As in the proof of 3.3(a) we see that for any character sheaf $A'$ not isomorphic to $A_i$ we have
$\c_i(A')=0$. From 2.4(b) we see that $\sum_j(-1)^j\c_i({}^pH^{j+2i}K_w)=1$ (we use that $\ce_i$ has rank $1$).
Hence $\sum_j(-1)^j(A_i:{}^pH^jK_w)\c_i(A_i)=1$ that is $(-1)^{\D+i}(A_i:{}^pH^{\D+i}K_w)\c_i(A_i)=1$.
Since $\c_i(A_i)=\pm1$ it follows that $(A_i:{}^pH^{\D+i}K_w)=1$ proving the first assertion of (b). The
remaining assertions of (b) follow from 2.3(a) applied with $y=w$ and with $A$ equal to $A_i$ or $A_{i'}$. 

\subhead 3.7\endsubhead 
In the setup of 3.6 we show:

(a) {\it for any $j$, $\op_k\ch^k({}^pH^jK_w)|_{\g_{i'}}$ does not contain $\ce_{i'}$ as a direct summand.}
\nl
Assume first that $j\ne\D+i$. It is enough to show that for any character sheaf $X$ such that
$(X:{}^pH^jK_w)\ne0$, $\op_k\ch^kX_{\g_{i'}}$ does not contain $\ce_{i'}$ as a direct summand. If $X$ is 
noncuspidal this follows from 1.6(a); if $X$ is cuspidal, it must be different from $A_i$ or $A_{i'}$ (see 
3.6(a)) and the result follows from the cleanness of cuspidal character sheaves other than $A_i$.

Assume now that for some $k$, $\ch^k({}^pH^{\D+i}K_w)|_{\g_{i'}}$ contains $\ce_{i'}$ as a direct summand. This, 
and the previous paragraph, imply that for some $k$, $\ch^k(K_w)|_{\g_{i'}}$ contains $\ce_{i'}$ as a direct 
summand. In particular $K_w|_{\g_{i'}}\ne0$ so that $\p_w\i(\g_{i'})\ne\em$. Using 3.1(a) we deduce that 
$\ul(w)\ge i'$ that is, $i\ge i'$. This contradiction proves (a).

\subhead 3.8\endsubhead  
We preserve the setup of 3.6. We have ${}^pH^{\D+i}K_w=Z+Z'$ since 

${}^pH^{\D+i}K_w/(Z+Z')$
\nl
 is both $A_i$-isotypic and $A_{i'}$-isotypic.
(It is a quotient of ${}^pH^{\D+i}K_w/Z$ which is $A_i$-isotypic and 
a quotient of ${}^pH^{\D+i}K_w/Z'$ which is $A_{i'}$-isotypic.)
As in the proof of 3.7(a) we see that all composition factors $X$ of $Z\cap Z'$ 
(which are necessarily not isomorphic to $A_i$ or $A_{i'}$) satisfy the condition that $\op_k\ch^k(X)|_{\g_{i'}}$
does not contain $\ce_{i'}$ as a direct summand. It follows that $\op_k\ch^k(Z\cap Z')_{\g_{i'}}$ does not 
contain $\ce_{i'}$ as a direct summand. Using this and 3.7(a) we deduce that 
$\op_k\ch^k({}^pH^{\D+i}K_w/(Z\cap Z'))_{\g_{i'}}$ does not contain $\ce_{i'}$ as a direct summand. 
Since ${}^pH^{\D+i}K_w=Z+Z'$, the natural map 
${}^pH^{\D+i}K_w/(Z\cap Z')@>>>({}^pH^{\D+i}K_w/Z)\op({}^pH^{\D+i}K_w/Z')$ is an isomorphism. It follows that

(a) $\op_k\ch^k({}^pH^{\D+i}K_w/Z)_{\g_{i'}}$ does not contain $\ce_{i'}$ as a direct summand;

(b) $\op_k\ch^k({}^pH^{\D+i}K_w/Z')_{\g_{i'}}$ does not contain $\ce_{i'}$ as a direct summand.
\nl
Since ${}^pH^{\D+i}K_w/Z\cong A_i$, we see that (a) (together with 1.5(a)) proves the cleanness of $A_i$ thus 
completing the proof of Theorem 1.2. 
Since ${}^pH^{\D+i}K_w/Z'$ is a direct sum of copies of $A_{i'}$ and 
$\op_k\ch^k(A_{i'})_{\g_{i'}}=\ce_{i'}$ we see that (b) implies ${}^pH^{\D+i}K_w/Z'=0$.
This, together with 3.6(a),(b) implies that

(c) {\it $(A_{i'}:{}^pH^jK_w)=0$ for any $j$.}

\subhead 3.9\endsubhead  
In view of the cleanness of $G$, we can restate 3.2(a) in a stronger form:

(a) {\it Let $i\ne i'$ in $I$. Then $\op_j\ch^jA_{i'}|_{\g_i}$ does not contain $\ce_i$ as a direct summand.}
\nl
Using this the proof of 3.3(a) applies in greater generality and yields the following result.

(b) {\it Assume that $i\in I$, $y\in\WW$, $\ul(y)<i$. Then $(A_i:\op_j{}^pH^jK_y)=0$. If, in addition, $i=10$ for
type $E_8$ and $i=4$ for type $F_4$ then $(A'_i:\op_j{}^pH^jK_y)=0$.}
\nl
From (b), 2.3(a) and 2.4(b) we deduce as in 3.6 the following result for any $i\in I$:

(c) {\it If $j\ne\D+i$ then $(A_i:{}^pH^jK_{w_i})=0$; if $j=\D+i$ then $(A_i:{}^pH^jK_{w_i})=1$ and there exists 
a unique subobject $Z$ of ${}^pH^jK_{w_i}$ such that $(A_i:Z)=0$ and ${}^pH^jK_{w_i}/Z\cong A_i$.}
\nl
The same result holds for $i=10$ (type $E_8$) and $i=4$ (type $F_4$) if $A_i$ is replaced by $A'_i$.

\subhead 3.10\endsubhead  
Note that, once Theorem 1.2 is known, the parity property \lb  \cite{\CS, III, (15.13.1)} can be established for a
reductive group in any characteristic as in \cite{\CS}. (Incidentally, note that 3.9(c) establishes the parity 
property for the character sheaves $A_i$ for $p=2$, type $E_8$ or $F_4$.) Using this we see that essentially the
same proof as in \cite{\CS} establishes \cite{\CS, V, Theorems 23.1, 24.4, 25.2, 25.6} (but not 
\cite{\CS, V, Theorem 24.8}) for a reductive group in any characteristic.

\enddocument